\documentclass{amsart}
\usepackage{amssymb}
\newtheorem{thm}{Theorem}
\newtheorem{lem}{Lemma}
\newtheorem{prop}{Proposition}
\newtheorem{cor}{Corollary}

\newtheorem{rem}{Remark}
\newtheorem{df}{Definition}

\begin{document}

\bibliographystyle{plain}

\title[Milo\v s Arsenovi\'c and Romi F. Shamoyan]{Trace theorems in harmonic function spaces on
$(\mathbb R^{n+1}_+)^m$, multipliers theorems and related problems}

\author[]{Milo\v s Arsenovi\' c$\dagger$}
\author[]{Romi F. Shamoyan}

\address{Department of mathematics, University of Belgrade, Studentski Trg 16, 11000 Belgrade, Serbia}
\email{\rm arsenovic@matf.bg.ac.rs}

\address{Bryansk University, Bryansk Russia}
\email{\rm rshamoyan@yahoo.com}

\thanks{$\dagger$ Supported by Ministry of Science, Serbia, project OI174017}

\date{}

\begin{abstract}
We introduce and study properties of certain new harmonic function spaces in products of upper half spaces. Norm estimates
for the so called expanded Bergman projection are obtained. Sharp theorems on multipliers acting on certain Sobolev type spaces of harmonic functions on the unit ball are obtained.
\end{abstract}

\maketitle

\footnotetext[1]{Mathematics Subject Classification 2010 Primary 42B15, Secondary 46E35.  Key words
and Phrases: harmonic functions, integral operators, traces, embedding theorems, unit ball, upper half space.}

\section{Introduction, preliminaries and auxiliary results}

The main goal of this paper is to introduce and study properties of certain new harmonic function spaces on the poly upper half space $(\mathbb R^{n+1}_+)^m$ and to solve trace problems for such spaces.
Solutions of trace problems in various functional spaces in complex function theory are based on estimates of Bergman type integral operators in various domains in $\mathbb C^n$, see for example \cite{DS}, \cite{JPS}, \cite{SM1}, \cite{SM3}, and references therein. In harmonic function spaces we used the same idea in \cite{AS4}. The second section of this paper provides some new estimates for such integral operators in the upper half space. Generally speaking, trace operator is acting on functions $f(z_1,...,z_m)$ defined on a product $\Omega^m$ of domains in $\mathbb R^k$, that is $z_j \in \Omega \subset
\mathbb R^k$, $1 \leq j \leq m$. However, when such a function $f$ is a product of $m$ functions  $f_1, \ldots, f_m$ we start to  deal with multi functional spaces. In the third section we provide a sharp embedding theorem for such spaces.

In the last section we give characterizations of the spaces of multipliers acting from Sobolev type mixed norm spaces of
harmonic functions on the unit ball into various spaces of harmonic functions on the unit ball.

We set $\mathbb H = \{(x, t) : x \in \mathbb R^n, t > 0 \} \subset \mathbb R^{n+1}$. For $z = (x, t) \in \mathbb H$ we set $\overline z = (x, -t)$. We denote the points in $\mathbb H$ usually by $z = (x, t)$ or $w = (y, s)$. The Lebegue measure is denoted by $dm(z) = dz = dx dt$ or $dm(w) = dw = dy ds$. We also use measures
$dm_\lambda(z) = t^\lambda dxdt$, $\lambda \in \mathbb R$.

We use common convention regarding constants: letter $C$ denotes a constant which can change its value from one occurrence
to the next one. Given two positive quantities $A$ and $B$, we write $A \asymp B$ if there are two constants $c, C > 0$ such that $cA \leq B \leq CA$.

The space of all harmonic functions in a domain $\Omega$ is denoted by $h(\Omega)$.
Weighted harmonic Bergman spaces on $\mathbb H$ are defined, for $0<p<\infty$ and $\lambda > -1$, as usual:
$$A^p_\lambda = A^p_\lambda(\mathbb H) = \left\{ f \in h(\mathbb H) : \| f \|_{A^p_\lambda} =
\left( \int_{\mathbb H} |f(z)|^p dm_\lambda(z) \right)^{1/p} < \infty \right\}.$$

For $\vec{\alpha} = (\alpha_1, \ldots, \alpha_m) \in \mathbb R^m$ we have a product measure $dm_{\vec{\alpha}}$ on 
$\mathbb H^m$ defined by $dm_{\vec{\alpha}}(z_1, \ldots, z_m) = dm_{\alpha_1}(z_1) \ldots dm_{\alpha_m}(z_m)$ and we
set $L^p_{\vec{\alpha}} = L^p(\mathbb H^m, dm_{\vec{\alpha}})$, $0<p<\infty$, and $A^p_{\vec{\alpha}} = L^p_{\vec{\alpha}} 
\cap h(\mathbb H^m)$. We denote by  $\tilde A^p_{\vec{\alpha}}$ the subspace of $A^p_{\vec{\alpha}}$ consisting of functions which are harmonic in each of the variables $z_1, \ldots, z_m$ separately.



We denote by $\mathbb B$ the open unit ball in $\mathbb R^n$ and by $\mathbb S = \partial \mathbb B$ the unit sphere in
$\mathbb R^n$. We denote polar coordinates in $\mathbb B$ by $x = rx'$, or $y = \rho y'$, where $x', y' \in \mathbb S$ and
$r = |x|$, $\rho = |y|$. Accordingly, the surface measure on $\mathbb S$ is denoted by $dx'$ or $dy'$.

Using multi index notation we set, for a function $f \in C^N(\Omega)$ and $N \in \mathbb N$:
\begin{equation*}
|\nabla^N f(x)| = \sqrt{ \sum_{|\gamma| = N} |D^\gamma f(x)|^2}, \qquad x \in \Omega.
\end{equation*}
For $0<p<\infty$, $0 \leq r < 1$ and $f \in C(\mathbb B)$ we set
$$M_p(f, r) = \left( \int_{\mathbb S} |f(rx')|^p dx' \right)^{1/p},$$
with the usual modification to cover the case $p = \infty$. For $0 < p \leq \infty$, $0 < q < \infty$, $\alpha > 0$
and $f \in C(\mathbb B)$ we consider mixed (quasi)-norms $\| f \|_{p,q,\alpha}$ defined by
\begin{equation}\label{qpnorm}
\| f \|_{p,q,\alpha} = \left( \int_0^1 M_p(f, r)^q (1-r^2)^{\alpha q - 1} r^{n-1}dr \right)^{1/q},
\end{equation}
again with the usual modification to cover the case $q = \infty$, and the corresponding spaces of harmonic functions
$$B^{p,q}_\alpha(\mathbb B) = B^{p,q}_\alpha = \{ f \in h(\mathbb B) : \| f \|_{p,q,\alpha} < \infty \}.$$
For details on these spaces we refer to \cite{DS}, Chapter 7. Also, for $N \in \mathbb N$ we have (quasi) norms
$$\| f \|_{D_NB^{p,q}_\alpha} = |f(0)| + \| \nabla^N f \|_{p,q,\alpha}, \qquad f \in C^N(\mathbb B),$$
and the corresponding spaces of harmonic functions:
\begin{equation*}
D_NB^{p,q}_\alpha = \{ f \in h(\mathbb B) : \| f \|_{D_NB^{p,q}_\alpha} < \infty \}.
\end{equation*}
We note that $A^p_\alpha = B^{p,p}_{\frac{\alpha + 1}{p}}$, therefore we have also spaces $D_NA^p_\alpha$. All of the above spaces are complete metric spaces, $D_NB^{p,q}_\alpha$ is a Banach space for $\min (p,q) \geq 1$ and $D_NA^p_\alpha$ is a Banach space for $p \geq 1$.



We also consider harmonic Triebel-Lizorkin spaces on the unit ball in $\mathbb R^n$, these were introduced in \cite{AS3} where embedding and multiplier results on these spaces can be found.
\begin{df}
Let $0 < p,q < \infty$ and $\alpha > 0$. The harmonic Triebel-Lizorkin space $F^{p,q}_\alpha (\mathbb B) =
F^{p,q}_\alpha$ consists of all functions $f \in h(\mathbb B)$ such that
\begin{equation}\label{tl}
\| f \|_{F^{p,q}_\alpha} = \left( \int_{\mathbb S} \left( \int_0^1 |f(rx')|^p (1-r)^{\alpha p -1} dr \right)^{q/p}
dx' \right)^{1/q} < \infty.
\end{equation}
\end{df}
These spaces are complete metric spaces, for $\min(p,q) \geq 1 $ they are Banach spaces.

Harmonic function spaces in the upper half spaces were studied recently in \cite{KK}, \cite{KY}, \cite{NY}, \cite{RY}.

\begin{df}
For a function $f : \mathbb H^m \to \mathbb C$ we define ${\rm Tr} f : \mathbb H \to \mathbb C$ by
${\rm Tr} f(z) = f(z,\ldots, z)$.

Let $X \subset h(\mathbb H^m)$. The trace of $X$ is ${\rm Trace}\; X = \{ {\rm Tr}\, f : f \in X \}$.
\end{df}

We denote the Poisson kernel for $\mathbb H$ by $P(x, t)$, i.e.
$$P(x, t) = c_n \frac{t}{(|x|^2 + t^2)^{\frac{n+1}{2}}}, \qquad x \in \mathbb R^n, \quad t > 0.$$
For $k \in \mathbb N_0$ a Bergman kernel $Q_k(z, w)$, where $z = (x, t) \in \mathbb H$ and $w = (y, s) \in \mathbb H$, is defined by
$$Q_k(z, w) = \frac{(-2)^{k+1}}{k!} \frac{\partial^{k+1}}{\partial t^{k+1}} P(x-y, t+s).$$
We need the following result from \cite{DS} which justifies terminology.

\begin{thm}\label{brthm}
Let $0<p<\infty$ and $\alpha > -1$. If $0<p\leq 1$ and $k \geq \frac{\alpha +n +1}{p} - (n+1)$ or $1\leq p < \infty$
and $k> \frac{\alpha +1}{p} - 1$, then
\begin{equation}\label{brep}
f(z) = \int_{\mathbb H} f(w) Q_k(z, w) s^k dy ds, \qquad f \in A^p_\alpha,\quad z \in \mathbb H.
\end{equation}
\end{thm}

The following elementary estimate of this kernel is contained in \cite{DS}:
\begin{equation}\label{estq}
|Q_k(z, w)| \leq C |z - \overline w|^{-(k+n+1)}, \quad z = (x, t), \quad w = (y, s) \in \mathbb H.
\end{equation}

Most of the results in the next two sections rely on the following three lemmas.

\begin{lem}[\cite{St}]\label{LemmaA}
There exists a collection $\{ \Delta_k \}_{k=1}^\infty$ of closed cubes in $\mathbb H$ with sides parallel to coordinate axes such that

$1^o$. $\cup_{k=1}^\infty \Delta_k = \mathbb H$ and ${\rm diam} \Delta_k \asymp {\rm dist} (\Delta_k,
\partial \mathbb H)$.

$2^o$.  The interiors of the cubes $\Delta_k$ are pairwise disjoint.

$3^o$.  If $\Delta_k^\ast$ is a cube with the same center as $\Delta_k$, but enlarged 5/4 times, then the
collection $\{ \Delta_k^\ast \}_{k=1}^\infty$ forms a finitely overlapping covering of $R^{n+1}_+$, i.e. there is
a constant $C = C_n$ such that $\sum_k \chi_{\Delta_k^\ast} \leq C$.
\end{lem}

\begin{lem}[\cite{D1}]\label{LemmaB}
Let $\Delta_k$ and $\Delta_k^\ast$ be the cubes from the previous lemma and let $(x_k, t_k)$ be the center of $\Delta_k$. Assume $f$ is subharmonic in $\mathbb H$. Then, for $0<p<\infty$ and $\alpha > 0$, we have
\begin{equation}
t_k^{\alpha p - 1} \max_{\Delta_k} |f|^p \leq \frac{C}{|\Delta_k^\ast|} \int_{\Delta_k^\ast} t^{\alpha p - 1}
|f(x,t)|^p dx dt, \qquad k \geq 1.
\end{equation}
\end{lem}

\begin{lem}[\cite{St}]\label{LemmaC}
Let $\Delta_k$ and $\Delta_k^\ast$ are as in the previous lemma, let $\zeta_k = (\xi_k, \eta_k)$ be the center of the cube $\Delta_k$. Then we have:
\begin{equation}\label{mlam}
 m_\lambda(\Delta_k) \asymp \eta_k^{n+1+\lambda} \asymp m_\lambda (\Delta_k^\ast), \qquad \lambda \in
\mathbb R,
\end{equation}
\begin{equation}\label{dis}
|\overline w - z | \asymp | \overline \zeta_k - z|, \qquad w \in \Delta_k^\ast, \quad z \in \mathbb H,
\end{equation}
\begin{equation}\label{height}
t \asymp \eta_k, \qquad (x,t) \in \Delta_k^\ast.
\end{equation}
\end{lem}


\begin{lem}[\cite{KY}]\label{omit}
If $\alpha > -1$ and $n+\alpha<2\gamma-1$, then
\begin{equation}
\int_{\mathbb H} \frac{t^\alpha dz}{|z - \overline w|^{2\gamma}} \leq C s^{\alpha + n + 1 - 2\gamma}, \qquad
w = (y, s) \in \mathbb H.
\end{equation}
\end{lem}

For $w = (y, s) \in \mathbb H$ we set $Q_w$ to be the cube, with sides parallel to the coordinate axis, centered at $w$ with
side length equal to $s$.

\section{Expanded Bergman projections and related operators}

In this section we provide new estimates for certain new integral operators closely connected with trace problem.

For any two $m$-tuples ($m \geq 1$) of reals $\vec{a} = (a_1, \ldots, a_m)$ and $\vec{b} = (b_1, \ldots, b_m)$ we define an
integral operator
\begin{equation}
(S_{\vec{a},\vec{b}}f)(z_1, \ldots, z_m) = \prod_{j=1}^m t_j^{a_j} \int_{\mathbb H}
\frac{f(w) s^{-n-1+\sum_{j=1}^m b_j}}{\prod_{j=1}^m |z_j - \overline w|^{a_j + b_j}}dw, \quad z_j = (x_j, t_j) \in \mathbb H.
\end{equation}
This operator can be called an expanded Bergman projection in the upper half space, it is well defined for $z_1, \ldots,
z_m \in \mathbb H$ and $f(w) \in L^1(\mathbb H, s^{-n-1-\sum_{j=1}^m b_j})$. A unit ball analogue of this operator was used in \cite{SM1}, see also \cite{Zh}.
We can write this operator in the following form:
$$(S_{\vec{a}, \vec{b}}f) (z_1, \ldots, z_m) = (N_{\vec{a}, \vec{b}}f) (z_1, \ldots, z_m) \prod_{j=1}^m t_j^{a_j},$$
where
$$N_{\vec{a}, \vec{b}} f (z_1, \ldots, z_m) = \int_{\mathbb H}
\frac{f(w) s^{-n-1+\sum_{j=1}^m b_j}}{\prod_{j=1}^m |z_j - \overline w|^{a_j + b_j}}dw.$$

We also consider related integral operators $S_{a,b}^k$, where $a>0$, $b>-1$ defined by
\begin{equation}
S_{a,b}^k f(z) = t^a \int_{\Delta_k} \frac{s^b f(w) dw}{|z-\overline w|^{n+1+a+b}},
\qquad z = (x, t) \in \mathbb H, \quad k \geq 1,
\end{equation}
and we set
\begin{equation}
\tilde S_{a,b} f(z) = S_{a,b}^k f(z), \qquad z \in \Delta_k.
\end{equation}

Analogous operators acting on analytic functions in the unit ball in $\mathbb C^n$ appeared in \cite{LS2}.

We need the following definition, generalizing the concept of Muckenhoupt weight to the upper half space. Analogous weights in the unit ball in $\mathbb C^n$ were considered in \cite{LS2} where a result analogous to Theorem \ref{sabbdd} was proven.

\begin{df}
Let $1<p<\infty$ and let $1/p + 1/q = 1$. A positive locally integrable function $V$ on $\mathbb H$ belongs to the $MH(p)$ class if
\begin{equation}
\sup_{w \in \mathbb H} \left( \frac{1}{|Q_w|} \int_{Q_w} V(z) dz \right) \left( \frac{1}{|Q_w|} \int_{Q_w} V(z)^{-q/p} dz
\right)^{p/q} < \infty.
\end{equation}
\end{df}

We remark here that an equivalent definition arises if in the above supremum we replace the family of cubes $Q_w$, $w \in \mathbb H$ with the family $\Delta_k$, $k \geq 1$. This easily follows from the fact that there is a constant $N= N_n$ such that each cube $Q_w$ can be covered by at most $N_n$ cubes from the family $\Delta_k$, and the selected cubes have sizes comparable to the size of $Q_w$.

Note that $V(z) = t^\alpha$ is in $MH(p)$ for all $1<p<\infty$ and all real $\alpha$.

\begin{thm}\label{sabbdd}
Let $0<\sigma<\infty$, $1 < p < \infty$ and $V \in MH(p)$. Then for every $f \in L^p_{\rm loc}(\mathbb H)$ we have
\begin{equation}
\sum_{k=1}^\infty \left( \int_{\Delta_k} |\tilde S_{a,b}^kf(z)|^p V(z) dm(z) \right)^{\sigma/p} \leq C \sum_{k=1}^\infty
\left(\int_{\Delta_k} |f(z)|^p V(z) dm(z) \right)^{\sigma/p}.
\end{equation}
\end{thm}

{\it Proof.} Let $q$ be the exponent conjugate to $p$. Let us fix $k \geq 1$. We have, using Lemma \ref{LemmaC} and Holder's inequality:
\begin{align*}
\phantom{\leq} & \int_{\Delta_k} |\tilde S_{a,b}^k f(z)|^p V(z) dm(z)
\leq  \int_{\Delta_k} \left( t^a \int_{\Delta_k} \frac{s^b |f(w)| dw}{|z-\overline w|^{n+1+a+b}} \right)^p V(z) dz \\
\leq & C \eta_k^{-p(n+1)} \int_{\Delta_k} V(z) dz \left( \int_{\Delta_k} |f(w)| V(w)^{1/p} V(w)^{-1/p} dw \right)^p \\
\leq & C \eta_k^{-p(n+1)} \int_{\Delta_k} V(z) dz \left( \int_{\Delta_k} V(w)^{-q/p} dw \right)^{p/q}
\int_{\Delta_k} |f(w)|^p V(w) dw \\
\leq & C \int_{\Delta_k} |f(w)|^p V(w) dw ,
\end{align*}
and this clearly proves the theorem. $\Box$


The following proposition is analogous to Proposition 1 from \cite{SM1}, the proof we present below follows the same pattern as the one provided in \cite{SM1} for the case of the unit ball in $\mathbb C^n$.

\begin{prop}\label{Pr3}
Let $1<p<\infty$, $a, b \in \mathbb R^m$ and $s_1, \ldots, s_m > -1$ satisfy $pa_j > -1-s_j$ and $p(mb_j -n) >
(m-1)(n+1) + ms_j + 1$ for $j=1,\ldots, m$. Set $\lambda = (m-1)(n+1) + \sum_{j=1}^m s_j$. Then there is a constant $C>0$ such that
\begin{equation}\label{eqPr2}
\int_{\mathbb H} \cdots \int_{\mathbb H} |(S_{\vec{a},\vec{b}}f)(z_1,\ldots,z_m)|^p dm_{s_1}(z_1) \ldots
dm_{s_m}(z_m) \leq C \| f \|^p_{L^p(\mathbb H, dm_\lambda)}
\end{equation}
for every $f \in L^p(\mathbb H, dm_\lambda)$.
\end{prop}

{\it Proof.} Let $1/p + 1/q = 1$. Choose $\gamma > 0$ such that
\begin{equation*}
p\gamma < p(mb_j - n) - (m-1)(n+1) - ms_j - 1, \qquad j = 1, \ldots, m.
\end{equation*}
Set $\alpha = \frac{1}{m} ( \gamma - \frac{1}{q})$ and choose $\beta$ such that $\beta + m\alpha = -n - 1 + \sum_{j=1}^m b_j$,
i.e. $\beta = -n -1 + \sum_{j=1}^m b_j - \gamma + \frac{1}{q}$. Since $pa_j + s_j + 1 > 0$ we can choose, for each $j = 1,
\ldots, m$, $e_j$ such that
\begin{equation*}
\frac{n+1}{mq} + \alpha < e_j < \frac{n+1}{mq} + \alpha + \frac{pa_j + s_j + 1}{p}.
\end{equation*}
Finally, set $d_j = a_j + b_j - e_j$.

After these preparations, we choose $f \in L^p(\mathbb H, dm_\lambda)$ and obtain, using Holder inequality with system of
$m+1$ exponents $p, mq, \ldots, mq$:
\begin{align*}
|N_{\vec{a}, \vec{b}} f (z_1, \ldots, z_m)| & \leq \int_{\mathbb H}
\frac{|f(w)| s^{-n-1+\sum_{j=1}^m b_j}}{\prod_{j=1}^m |z_j - \overline w|^{a_j + b_j}}dw \\
& = \int_{\mathbb H} \frac{|f(w)| s^\beta}{\prod_{j=1}^m |z_j - \overline w|^{d_j}} \prod_{j=1}^m \frac{s^\alpha}{|z_j-
\overline w|^{e_j}} dw \\
& \leq \left( \int_{\mathbb H} \frac{|f(w)|^p s^{p\beta}dw}{\prod_{j=1}^m |z_j - \overline w|^{pd_j}} \right)^{1/p}
\prod_{j=1}^m \left( \int_{\mathbb H} \frac{s^{qm\alpha}}{|z_j-\overline w|^{qme_j}} dw \right)^{\frac{1}{qm}}\\
& \leq C \left( \int_{\mathbb H} \frac{|f(w)|^p s^{p\beta}dw}{\prod_{j=1}^m |z_j - \overline w|^{pd_j}} \right)^{1/p}
\prod_{j=1}^m t_j^{\alpha - e_j + \frac{n+1}{qm}},
\end{align*}
where, at the last step, we used Lemma \ref{omit}. Therefore we have
\begin{equation*}
|S_{\vec{a}, \vec{b}} f(z_1, \ldots, z_m)|^p \leq C \int_{\mathbb H} \frac{|f(w)|^p s^{p\beta}dw}{\prod_{j=1}^m
|z_j - \overline w|^{pd_j}} \prod_{j=1}^m t_j^{p(a_j+\alpha - e_j + \frac{n+1}{qm})}.
\end{equation*}
Hence, using Fubini's theorem and Lemma \ref{omit} we obtain
\begin{align*}
\phantom{\leq} & \int_{\mathbb H} \cdots \int_{\mathbb H} |(S_{\vec{a},\vec{b}}f)(z_1,\ldots,z_m)|^p dm_{s_1}(z_1) \ldots
dm_{s_m}(z_m)\\
\leq & C \int_{\mathbb H} |f(w)|^ps^{p\beta} \left( \int_{\mathbb H} \cdots \int_{\mathbb H} \prod_{j=1}^m
\frac{t_j^{s_j + p(a_j + \alpha -e_j + \frac{n+1}{qm})}}{|z_j - \overline w|^{pd_j}} dz_1 \ldots dz_m \right) dw\\
= & C \int_{\mathbb H} |f(w)|^ps^{p\beta} \left( \prod_{j=1}^m \int_{\mathbb H} \frac{t_j^{s_j + p(a_j + \alpha -e_j + \frac{n+1}{qm})}}{|z_j - \overline w|^{pd_j}} dz_j \right) dw = C \int_{\mathbb H} |f(w)|^p s^\theta dw,
\end{align*}
where we have, see Lemma \ref{omit},
\begin{align*}
\theta & = p \beta + \sum_{j=1}^m \left[ s_j + p(a_j + \alpha - e_j + (n+1)/qm) - pd_j + n + 1 \right]\\
& = \sum_{j=1}^m s_j + m(n+1) + \frac{p(n+1)}{q} + p \left( \beta + m\alpha + \sum_{j=1}^m (a_j - e_j - d_j)
\right)\\
& = \sum_{j=1}^m s_j + (n+1)(m + \frac{p}{q}) + p \left( -n-1 + \sum_{j=1}^m (b_j + a_j - e_j-d_j) \right) = \lambda
\end{align*}
and this ends the proof. $\Box$

Next we consider another class of integral operators, see \cite{SM2} for similar operators acting on analytic functions in
poly balls and for an analogue of Proposition \ref{aadmpr} below.

For any two $m$-tuples $\vec{a} = (a_1, \ldots, a_m)$ and $\vec{b} = (b_1, \ldots, b_m)$ of reals we set
\begin{equation*}\label{rop1}
(R_{\vec{a}, \vec{b}} g)(w) = s^{-m(n+1) + \sum_{j=1}^m b_j} \int_{\mathbb H} \cdots \int_{\mathbb H} g(z_1, \ldots, z_m)
\prod_{j=1}^m \frac{t_j^{a_j}}{|z_j - \overline w|^{a_j + b_j}} dz_1 \ldots dz_m,
\end{equation*}
where $w = (y, s) \in \mathbb H$ and $g \in L^1_{\vec{a}}$. Next, for $k \in \mathbb N_0$ we define an integral operator
\begin{equation*}\label{rop2}
(R_k g)(w) = \int_{\mathbb H} \cdots \int_{\mathbb H} g(z_1, \ldots, z_m)
\prod_{j=1}^m Q_k(z_j,w) dm_k(z_1) \ldots dm_k(z_m), \qquad w \in \mathbb H.
\end{equation*}

In fact, this operator is the trace operator on a suitable space. Indeed we have
\begin{equation*}
(R_k g)(w) = g(w), \qquad g \in \tilde A^p_{\vec{\alpha}}, \quad \alpha = (\alpha_1, \ldots, \alpha_m),
\end{equation*}
if $p$, $n$ and $\alpha_j$ satisfy conditions from Theorem \ref{brthm}.


The following proposition is well known in the case of analytic functions in the unit ball in $\mathbb C^n$, see \cite{Zh},
it was extended to analytic functions on poly balls in \cite{SM2} and here we deal with harmonic functions in the poly
half space.

\begin{prop}\label{aadmpr}
Let $1 \leq p < \infty$ and $\vec{a}, \vec{b}, \vec{\alpha} \in \mathbb R^m$. If $p>1$ we assume these parameters satisfy the following conditions:
$$q(a_j - \alpha_j) > -1 - \alpha_j, \qquad 1 \leq j \leq m,$$
$$q(m(b_j + \alpha_j)-n) > (m-1)(n+1) + m\alpha_j + 1, \qquad 1 \leq j \leq m,$$
where $q$ is the exponent conjugate to $p$. If $p=1$ we assume $m(\alpha_j + b_j)> n$ and $\alpha_j < a_j$ for
$j = 1, \ldots, m$. Set $\lambda = (m-1)(n+1) + \sum_{j=1}^m \alpha_j$. Then
\begin{equation}\label{Rbdd}
\| R_{\vec{a}, \vec{b}} g \|_{L^p(\mathbb H, dm_\lambda)} \leq C \| g \|_{L^p_{\vec{\alpha}}}, \qquad g \in L^p_{\vec{\alpha}}.
\end{equation}
\end{prop}

{\it Proof.} We follow the same method as in \cite{SM2}, adapted to our situation. Let us start with the case $p=1$. By Fubini's theorem we have
\begin{align}
\phantom{\leq} & \| R_{\vec{a}, \vec{b}} g \|_{L^1(\mathbb H, dm_\lambda)} = \int_{\mathbb H} |(R_{\vec{a}, \vec{b}}g)(w)| s^\lambda dw  \notag\\
\leq & \int_{\mathbb H} \cdots \int_{\mathbb H} |g(z_1, \ldots, z_m)|
\prod_{j=1}^m t_j^{a_j} \int_{\mathbb H} \frac{s^{-n-1 + \sum_{j=1}^m (\alpha_j+b_j)} dw}{\prod_{j=1}^m |z_j - \overline w
|^{a_j+b_j}} dz_1 \ldots dz_m. \label{combin}
\end{align}
Next we use Holder's inequality with $m$ functions and Lemma \ref{omit} to obtain:
\begin{align*}
\int_{\mathbb H} \frac{s^{-n-1 + \sum_{j=1}^m (\alpha_j+b_j)} dw}{\prod_{j=1}^m |z_j - \overline w|^{a_j+b_j}}
\leq & \prod_{j=1}^m \left( \int_{\mathbb H} \frac{s^{-n-1 + m(\alpha_j + b_j)} dw}{|z_j - \overline w|^{m(a_j+b_j)}}
\right)^{1/m}\\
\leq & C \prod_{j=1}^m t_j^{\alpha_j - a_j}.
\end{align*}
This estimate, combined with (\ref{combin}) settles the case $p=1$.

Next we assume $1<p<\infty$. Let $q$ be the exponent conjugate to $p$. Using identity
\begin{align*}
\phantom{=} & \int_{\mathbb H} (R_{\vec{a}, \vec{b}}g)(w) f(w) dm_\lambda (w)\\
=  & \int_{\mathbb H} \left( \int_{\mathbb H} \cdots \int_{\mathbb H} g(z_1, \ldots, z_m) \prod_{j=1}^m \frac{t_j^{a_j}}{|z_j - \overline w|^{a_j + b_j}} dz_1 \ldots dz_m \right) \\
\phantom{=} & f(w) s^{-n-1 + \sum_{j=1}^m (\alpha_j + b_j)}dw\\
= & \int_{\mathbb H} \cdots \int_{\mathbb H} g(z_1, \ldots, z_m) \left( \int_{\mathbb H} s^{-n-1+ \sum_{j=1}^m(\alpha_j + b_j)} \frac{\prod_{j=1}^m t_j^{a_j - \alpha_j} f(w)dw}{\prod_{j=1}^m |z_j - \overline w|^{a_j+b_j}} \right)\\
\phantom{=} & dm_{\alpha_1}(z_1) \ldots dm_{\alpha_m}(z_m),
\end{align*}
valid for, for example, continuous compactly supported $f \in L^q(dm_\lambda)$ and $g \in L^p_{\vec{\alpha}}$ we see that
the conjugate operator $R_{\vec{a}, \vec{b}}^\star$ is equal to $S_{\vec{a} - \vec{\alpha}, \vec{b} + \vec{\alpha}}$.
However, the last one is bounded from $L^q(dm_\lambda)$ to $L^q_{\vec{\alpha}}$ by Proposition \ref{Pr3} and therefore
$R_{\vec{a}, \vec{b}} : L^p_{\vec{\alpha}} \to L^p(dm_\lambda)$ is also bounded. $\Box$

Using (\ref{estq}) we see that $|R_k g(w)| \leq (R_{\vec{a}, \vec{b}} |g|)(w)$ where $a_j = k$ and $b_j = n+1$ for
$j = 1, \ldots, m$. This observation leads to the following corollary.

\begin{cor}
Let $k \in \mathbb N_0$, $1 \leq p < \infty$ and $\alpha_j > -1$ for $j = 1, \ldots, m$. If $p=1$ we assume $\alpha_j < k$
for $1 \leq j \leq m$, if $1<p<\infty$ we assume
$$q(k-\alpha_j) > -1 - \alpha_j, \qquad 1 \leq j \leq m,$$
$$q(m(n+1+\alpha_j)-n) > (m-1)(n+1) + m\alpha_j + 1, \qquad 1 \leq j \leq m.$$
Set $\lambda = (m-1)(n+1) + \sum_{j=1}^m \alpha_j$. Then the operator $S_k$ maps $L^p_{\vec{\alpha}}$ continuously into $L^p(\mathbb H, dm_\lambda)$
\end{cor}

\section{Trace theorems and embedding theorems for multi functional spaces of harmonic functions}

In this section we give an estimate of the $A^p_\lambda$-norm of trace, Theorem \ref{trest} below.  Theorem \ref{tr1T} is a sharp embedding result obtained using norm estimate of the operator $S_{\vec{a}, \vec{b}}$, while Theorem \ref{trised} is a sharp embedding theorem closely connecting trace operator and multi functional spaces. At the end of this section we consider Carleson type conditions adapted to multi functional setting for positive Borel measures on poly upper half spaces,
see Definition \ref{drcar} and Theorem \ref{Thrcar}.


\begin{lem}[\cite{AS4}]\label{Lemma4}
Let $0<p<\infty$ and $s_1, \ldots, s_m > -1$. Set $\lambda = (m-1)(n+1) + \sum_{j=1}^m s_j$. Then there is a constant
$C>0$ such that for all $f \in  h(\mathbb H^m)$ we have
\begin{equation}\label{eqL4}
\int_{\mathbb H} |{\rm Tr}\, f(z)|^p dm_\lambda (z) \leq C \int_{\mathbb H} \cdots \int_{\mathbb H}
|f(z_1, \ldots, z_m)|^p dm_{s_1}(z_1) \ldots dm_{s_m}(z_m).
\end{equation}
\end{lem}

A holomorphic version of the following theorem appeared in \cite{LS1}.
\begin{thm}\label{trest}
Let $\alpha > -1$, assume $f_i \in h(\mathbb H^t)$ for $i = 1, \ldots, m$. Let $0 < p_i, q_i < \infty$ satisfy $\sum_{i=1}^m
\frac{p_i}{q_i} = 1$ and assume $\beta_i = \frac{(n+1+\alpha)q_i}{tmp_i} - (n+1) > -1$ for $1 \leq i \leq m$. Then
\begin{align}
\phantom{\leq} & \int_{\mathbb H} |{\rm Tr}\, f_1(w)|^{p_1} \ldots |{\rm Tr}\, f_m(w)|^{p_m} s^\alpha dw \label{bezi} \\
\leq & C \prod_{i=1}^m \left( \int_{\mathbb H} \cdots \int_{\mathbb H} |f_i(w_1, \ldots, w_t)|^{q_i} \prod_{j=1}^t
s_j^{\beta_i} dw_1 \ldots dw_t \right)^{p_i/q_i}. \notag
\end{align}
\end{thm}

{\it Proof.} Let us denote the integral appearing in (\ref{bezi}) by $I$. Using Lemma \ref{LemmaA} and Lemma \ref{LemmaC} we obtain
\begin{align}
I & = \sum_{k=1}^\infty \int_{\Delta_k} \prod_{i=1}^m |f_i(w, \ldots, w)|^{p_i} s^\alpha dw \leq C \sum_{k=1}^\infty
\eta_k^{n+1+\alpha} \sup_{w \in \Delta_k} \prod_{i=1}^m |f_i(w, \ldots, w)|^{p_i} \notag \\
& \leq C \sum_{k=1}^\infty \eta_k^{n+1+\alpha} \prod_{i=1}^m \sup_{w \in \Delta_k} |f_i(w, \ldots, w)|^{p_i} \notag\\
& \leq C \sum_{k=1}^\infty \eta_k^{n+1+\alpha} \prod_{i=1}^m \sup_{w_1, \ldots, w_t \in \Delta_k} |f_i(w_1, \ldots,
w_t)|^{p_i} \label{single}
\\
& \leq C \sum_{k_1 = 1}^\infty \cdots \sum_{k_t = 1}^\infty \left( \sup_{w_j \in \Delta_{k_j}} |f_1(w_1, \ldots, w_t)|^{p_1}
\cdots  \sup_{w_j \in \Delta_{k_j}} |f_m(w_1, \ldots, w_t)|^{p_m} \right) \notag \\
& \phantom{\leq} \times \eta_{k_1}^{\frac{n+1+\alpha}{t}}  \cdots  \eta_{k_t}^{\frac{n+1+\alpha}{t}} \label{double}.
\end{align}
The last inequality follows from the fact that for $k_1 = \cdots = k_m = k$ expression in (\ref{double}) reduces to
(\ref{single}). Next we apply generalized Holder's inequality with exponents $q_i/p_i$, $1 \leq i \leq m$,
to the last multiple sum and obtain
\begin{align}
I & \leq C \left( \sum_{k_1, \ldots, k_t = 1}^\infty \sup_{w_j \in \Delta_{k_j}} |f_1(w_1, \ldots, w_t)|^{q_1}
(\eta_{k_1} \ldots \eta_{k_t})^{\frac{(n+1+\alpha)q_1}{mtp_1}} \right)^{p_1/q_1}\times \ldots \notag\\
& \phantom{\leq} \times  \left( \sum_{k_1, \ldots, k_t = 1}^\infty \sup_{w_j \in \Delta_{k_j}} |f_m(w_1, \ldots, w_t)|^{q_m}
(\eta_{k_1}\ldots \eta_{k_t})^{\frac{(n+1+\alpha)q_m}{mtp_m}} \right)^{p_m/q_m}. \label{incom} \end{align}
Since
\begin{equation}
\sup_{w_j \in \Delta_{k_j}} |f_j(w_1, \ldots, w_t)|^{q_j} \leq C \int_{w_1 \in \Delta_{k_1}^\ast} \ldots
\int_{w_t \in \Delta_{k_t}^\ast} |f_j|^{q_j} dw_1 \ldots dw_t \prod_{i=1}^t \eta_{k_i}^{-n-1} \end{equation}
we obtain, using finite overlapping property of $\Delta_k^\ast$:

\begin{align*}
I_j & = \sum_{k_1, \ldots, k_t = 1}^\infty \sup_{w_j \in \Delta_{k_j}} |f_j(w_1, \ldots, w_t)|^{q_j}
(\eta_{k_1}\ldots \eta_{k_t})^{\frac{(n+1+\alpha)q_j}{mtp_j}}\\
& \leq C  \sum_{k_1, \ldots, k_t = 1}^\infty
 \int_{w_1 \in \Delta_{k_1}^\ast} \ldots
\int_{w_t \in \Delta_{k_t}^\ast} |f_j|^{q_j} dw_1 \ldots dw_t \prod_{i=1}^t \eta_{k_i}^{\beta_j}\\
& \leq C  \sum_{k_1, \ldots, k_t = 1}^\infty \int_{w_1 \in \Delta_{k_1}^\ast} \ldots \int_{w_t \in \Delta_{k_t}^\ast}
|f_j|^{q_j} s_1^{\beta_j} \ldots s_t^{\beta_j} dw_1 \ldots dw_t\\
& \leq C \int_{\mathbb H} \cdots \int_{\mathbb H} |f_j(w_1, \ldots, w_t)|^{q_j} (s_1 \ldots s_t)^{\beta_j} dw_1 \ldots
dw_t.
\end{align*}
This, in combination with (\ref{incom}), suffices to establish needed estimate. $\Box$

The theorem below was announced, without proof, in \cite{SM1} for $0<p<\infty$. A proof for the case $0<p\leq 1$ was given in \cite{AS4}, here we settle the remaining case $1<p<\infty$.

\begin{thm}\label{tr1T}
Let $1<p<\infty$, $s_1, \ldots, s_m > -1$ and set $\lambda = (m-1)(n+1) + \sum_{j=1}^m s_j$. Then
\begin{equation}\label{eqtr1}
A^p_\lambda \subset {\rm Trace}\; \tilde A^p_{\overrightarrow{s}} \subset {\rm Trace}\; A^p_{\overrightarrow{s}} \subset
L^p(\mathbb R^{n+1}_+, dm_\lambda).
\end{equation}
In particular, if $f \in A^p_{\overrightarrow{s}}$ and if ${\rm Tr}\; f$ is harmonic, then ${\rm Tr}\; f \in A^p_\lambda$.
\end{thm}

{\it Proof.} The second inclusion in (\ref{eqtr1}) is trivial while the third one follows from Lemma \ref{Lemma4}. Let us
prove the first inclusion, we fix $g \in A^p_\lambda$. Let us choose $k \in \mathbb N_0$ such that $p(n+k+1) > (m-1)(n+1) +
ms_j + pn + 1$ for $1\leq j \leq m$ and set
\begin{equation*}
f(z_1, \ldots, z_m) = \int_{\mathbb H} Q_k(\frac{z_1 + \cdots + z_m}{m}, w) g(w) s^k dw,
\qquad z_1, \ldots, z_m \in \mathbb H.
\end{equation*}

We have, by Theorem \ref{brthm}, ${\rm Tr}\, f = g$. Since the kernel $Q_k(m^{-1}(z_1 + \cdots + z_m), w)$
is harmonic in each of the variables $z_1, \ldots, z_m$ it follows that the same is true for $f(z_1, \ldots, z_m)$. Using
estimate (\ref{estq}) and classical inequality between arithmetic and geometric mean we obtain
\begin{equation}
|f(z_1, \ldots, z_m)|  \leq C \int_{\mathbb H}
\frac{|g(w)|s^k dw}{\left| \frac{z_1 + \cdots + z_m}{m} - \overline w \right|^{k+n+1}}
\leq C \int_{\mathbb H} \frac{|g(w)|s^k dw}{ \prod_{j=1}^m |z_j - \overline w|^{\frac{n+k+1}{m}}}.
\end{equation}
Hence $|f(z_1, \ldots, z_m)| \leq C (S_{\vec{a}, \vec{b}}\,|g|)(z_1, \ldots, z_m)$ where $a_j = 0$ and $b_j =
(n+k+1)/m$ for $j = 1, \ldots, m$. Now an application of Proposition \ref{Pr3} completes the proof. $\Box$

\begin{lem}\label{eqset}
Let $0<q,\sigma<\infty$ and $\alpha > -1$. Then we have
\begin{equation*}
\sum_{k=1}^\infty \eta_k^{n+1} \left( \int_{\Delta_k} |f(z)|^\sigma dm_\alpha (z) \right)^{\frac{q}{\sigma}} \leq C
\int_{\mathbb H} \left( \int_{Q_w} |f(z)|^\sigma dm_\alpha (z) \right)^{\frac{q}{\sigma}} dw, \quad f \in h(\mathbb H).
\end{equation*}
\end{lem}

This is a special case of Lemma 6 from \cite{SA1}, in fact harmonicity of $f$ is not needed here.

\begin{lem}\label{elsum}
Let $0<q_i\leq p < \infty$ and let $x_{i,k} \geq 0$ for $1 \leq i \leq m$ and $k \geq 1$. Then:
\begin{equation*}
\left( \sum_{k=1}^\infty x_{1,k}^p x_{2,k}^p \ldots x_{m,k}^p \right)^{1/p} \leq \prod_{i=1}^m \left( x_{i,1}^{q_i} +
x_{i,2}^{q_i} + \cdots \right)^{1/q_i}.
\end{equation*}
\end{lem}

{\it Proof.} Since the $l^q$ norm of a sequence is a decreasing function of $q$ we can assume $q_i = p$ for all $i = 1, \ldots, m$. But in this special case our inequality is equivalent to
\begin{equation*}
 \sum_{k=1}^\infty x_{1,k}^p x_{2,k}^p \ldots x_{m,k}^p  \leq \prod_{i=1}^m \left( x_{i,1}^p +
x_{i,2}^p + \cdots \right)
\end{equation*}
and this is clearly true. $\Box$

\begin{thm}\label{trised}
Let $0<p<\infty$ and let, for $i=1,\ldots, m$, $0< q_i, \sigma_i, <\infty$, $\alpha_i > -1$. Assume $q_i \leq p$ for $i = 1, \ldots, m$. Let $\mu$ be a positive Borel measure on $\mathbb H$. Then the following two conditions are equivalent:

$1^o$. For any harmonic function $f(z_1, \ldots, z_m)$ on $\mathbb H^m$ that splits into a product of harmonic functions
$f_i(z_i) \in h(\mathbb H)$, i.e $f(z_1, \ldots, z_m) = \prod_{j=1}^m f_j(z_j)$ we have
\begin{equation}\label{trsteqn}
\left( \int_{\mathbb H} |{\rm Tr}\, f(z)|^p d \mu(z)\right)^{1/p} \leq C \prod_{i=1}^m \left( \int_{\mathbb H}
\left( \int_{Q_w} |f_i(z)|^{\sigma_i} dm_{\alpha_i}(z) \right)^{q_i/\sigma_i} dw \right)^{1/q_i}.
\end{equation}

$2^o$. The measure $\mu$ satisfies the following Carleson-type condition:
\begin{equation}\label{cacoqn}
\mu(\Delta_k) \leq C \eta_k^{p\sum_{i=1}^m \left(\frac{n+1+\alpha_i}{\sigma_i} + \frac{n+1}{q_i}\right)},
\qquad k \geq 1.
\end{equation}
\end{thm}

{\it Proof.} Set $\theta = p\sum_{i=1}^m \left(\frac{n+1+\alpha_i}{\sigma_i} + \frac{n+1}{q_i}\right)$ and $\theta_i =
p\left( \frac{n+1+\alpha_i}{\sigma_i} + \frac{n+1}{q_i}\right)$, $1 \leq i \leq m$. Using a covering argument one easily shows that for any positive measurable function $u : \mathbb H \to \mathbb R$ we have
\begin{equation}\label{gen}
\sum_{k=1}^\infty \eta_k^{n+1} \left( \int_{\Delta_k^\ast} u(z) dm_\alpha(z) \right)^\beta \asymp
\sum_{k=1}^\infty \eta_k^{n+1} \left( \int_{\Delta_k} u(z) dm_\alpha(z) \right)^\beta, \qquad \beta > 0.
\end{equation}

Let us assume
(\ref{cacoqn}) holds. Then we have, using Lemma \ref{LemmaA} and Lemma \ref{LemmaB}
\begin{align*}
\phantom{\leq} & \int_{\mathbb H} |{\rm Tr}\, f(z)|^p d \mu(z) = \sum_{k=1}^\infty \int_{\Delta_k}
|{\rm Tr}\, f(z)|^p d \mu(z)
\leq \sum_{k=1}^\infty \mu(\Delta_k) \max_{z \in \Delta_k}
|{\rm Tr}\, f(z)|^p\\
\leq & C \sum_{k=1}^\infty \eta_k^\theta \prod_{i=1}^m \left( \max_{z \in \Delta_k} |f_i(z)|^{\sigma_i} \right)^{p/\sigma_i}
=  C \sum_{k=1}^\infty \prod_{i=1}^m \eta_k^{\theta_i} \left( \max_{z \in \Delta_k} |f_i(z)|^{\sigma_i} \right)^{p/\sigma_i}\\
\leq & C \sum_{k=1}^\infty \prod_{i=1}^m \eta_k^{\theta_i} \left( \eta_k^{-n-1-\alpha_i}
\int_{\Delta_k^\ast} |f_i(z)|^{\sigma_i} dm_{\alpha_i}(z) \right)^{p/\sigma_i}\\
= & C \sum_{k=1}^\infty \prod_{i=1}^m \eta_k^{p\frac{n+1}{q_i}} \left( \int_{\Delta_k^\ast} |f_i(z)|^{\sigma_i} dm_{\alpha_i}(z) \right)^{p/\sigma_i}\\
= & C \sum_{k=1}^\infty x_{1,k}^p x_{2,k}^p \ldots x_{m,k}^p, \quad x_{i,k} = \eta_k^{\frac{n+1}{q_i}}
\left( \int_{\Delta_k^\ast} |f_i(z)|^{\sigma_i} dm_{\alpha_i}(z) \right)^{1/\sigma_i}.
\end{align*}
Now an application of Lemma \ref{elsum} followed by (\ref{gen}) gives
\begin{align*}
\left( \int_{\mathbb H} |{\rm Tr}\, f(z)|^p d \mu(z) \right)^{1/p} & \leq C \prod_{i=1}^m \left( \sum_{k = 1}^\infty
\eta_k^{n+1} \left( \int_{\Delta_k^\ast} |f_i(z)|^{\sigma_i} dm_{\alpha_i}(z) \right)^{q_i/\sigma_i}\right)^{1/q_i} \\
& \leq C \prod_{i=1}^m \left( \sum_{k = 1}^\infty
\eta_k^{n+1} \left( \int_{\Delta_k} |f_i(z)|^{\sigma_i} dm_{\alpha_i}(z) \right)^{q_i/\sigma_i}\right)^{1/q_i},
\end{align*}
which, in view of Lemma \ref{eqset}, is sufficient to derive (\ref{trsteqn}).

We give an outline of proof of the reverse implication. Namely, one uses test functions
\begin{equation*}
f(z_1, \ldots, z_m) = \prod_{j=1}^m f_j(z_j), \qquad f_j(z) = f_{\theta_k, l}(z) = \frac{\partial^l}{\partial t^l}
\frac{1}{|z - \overline \theta_k |^{n-1}}, \quad 1 \leq j \leq m,
\end{equation*}
where $l$ is sufficiently large, $\theta_k$ is suitably chosen point near $\zeta_k$, see \cite{AS4} and the proof of
Theorem 4 from \cite{SA1} for these choices. The right hand side of (\ref{trsteqn}) is estimated with help of pointwise
estimates of $f_{w,l}$ from \cite{AS4} and Lemma \ref{omit}. We leave details to the interested reader. $\Box$







\begin{df}\label{drcar}
Let $\mu$ be a positive Borel measure on $\mathbb H^m$ and let $r_1, \ldots, r_m > 0$. We say $\mu$ is an
$(r_1, \ldots, r_m)$-Carleson measure if
\begin{equation}\label{rcarm}
\| \mu \|_{(r_1,\ldots,r_m)} = \sup_{w_1, \ldots, w_m \in \mathbb H}
\frac{\mu(Q_{w_1} \times \cdots \times Q_{w_m})}{s_1^{r_1} \ldots s_m^{r_m}} < \infty, \qquad w_j=(y_j,s_j).
\end{equation}
\end{df}

The following theorem, which is an analogue of Theorem 2 from \cite{SM1}, see also \cite{Zh}, gives an equivalent
description of $(r_1, \ldots, r_m)$-Carleson measures.

\begin{thm}\label{Thrcar}
Let $\mu$ be a positive Borel measure on $\mathbb H^m$. Assume $r_1, \ldots, r_m > n$ and $\tau_1, \ldots, \tau_m > 0$.
Then $\mu$ is an $(r_1, \ldots, r_m)$-Carleson measure
if and only if
\begin{equation}\label{carin}
\| \mu \|_{(r_1,\ldots,r_m)}^\ast = \sup_{w_1, \ldots, w_m \in \mathbb H} \int_{\mathbb H_{w_1}} \ldots
\int_{\mathbb H_{w_m}}
\prod_{j=1}^m \frac{s_j^{\tau_j}}{|z_j - \overline w_j|^{r_j+\tau_j}} d\mu(z_1, \ldots, z_m) < \infty,
\end{equation}
where $\mathbb H_{w_j} = \{ w\in \mathbb H : s \leq 3s_j \}$.
Moreover, $\| \mu \|_{(r_1,\ldots,r_m)} \asymp \| \mu \|_{(r_1,\ldots,r_m)}^\ast$.
\end{thm}

{\it Proof.} Assume (\ref{carin}) holds and choose $w_1, \ldots, w_m \in \mathbb H$. Then we have, using Lemma \ref{LemmaC},
\begin{align*}
\| \mu \|_{(r_1,\ldots,r_m)}^\ast & \geq \int_{Q_{w_1}} \cdots \int_{Q_{w_m}} \prod_{j=1}^m \frac{s_j^{\tau_j}}{|z_j -
\overline w_j|^{r_j + \tau_j}} d\mu(z_1, \ldots, z_m)\\
& \geq C \frac{\mu(Q_{w_1} \times \cdots Q_{w_m})}{s_1^{r_1} \ldots
s_m^{r_m}},
\end{align*}
which implies that $\mu$ is an $(r_1, \ldots, r_m)$-Carleson measure. Note that in this implication we did not use
conditions on the parameters.

Now we assume $\mu$ is an $(r_1, \ldots, r_m)$-Carleson measure. Let us, moreover, assume $m=1$. We choose $w = (y, s) \in
\mathbb H$, in order to simplify notation we assume $y = 0$. Let $\Gamma = \mathbb Z^n$ be the integer lattice in $\mathbb R^n$. We have a partition of $\mathbb H_w$ into layers $H_{k,s} = \{ z \in \mathbb H : 2^{-k}s \leq t < 3\cdot 2^{-k}s \}$,
$k \in \mathbb N_0$. Moreover, each layer $H_{k,s}$ is partitioned into congruent cubes $Q_{k, \xi_j}$ with centers
$\theta_{k,j} = (2^{-k}s \xi_j, 2^{-k}s)$, where $\xi_j \in \Gamma$. Since
\begin{equation}
|z-(0,-s)| \asymp \sqrt{s^2 + |2^{-k}s\xi_j|^2}, \qquad z \in Q_{k, \xi_j}
\end{equation}
we obtain
\begin{align*}
\int_{\mathbb H} \frac{s^\tau d\mu(z)}{|z-\overline w|^{r+\tau}} & = \sum_{k=0}^\infty \int_{\mathbb H_{k,s}} \frac{s^\tau d\mu(z)}{|z-\overline w|^{r+\tau}} = s^{\tau} \sum_{k = 0}^\infty \sum_{\xi_j \in \Gamma} \int_{Q_{k, \xi_j}}
\frac{d\mu(z)}{|z-\overline w|^{r+\tau}}\\
& \leq \| \mu \|_r s^{\tau} \sum_{k=0}^\infty \sum_{\xi_j \in \Gamma} \frac{(2^{-k}s)^r}{(s^2 + |2^{-k}s\xi_j|^2)^{\frac{r+\tau}{2}}}\\
& =  \| \mu \|_r \sum_{k=0}^\infty \sum_{\xi_j \in \Gamma}
\frac{(2^{-k})^r}{(1 + |2^{-k}\xi_j|^2)^{\frac{r+\tau}{2}}}\\
& \leq C \| \mu \|_r \sum_{k=0}^\infty 2^{-kr} \int_{\mathbb R^n} \frac{dx}{(1+|2^{-k}x|^2)^{r+\tau}} \\
& = C \| \mu \|_r \sum_{k=0}^\infty 2^{-kr} \int_0^\infty \frac{r^{n-1} dr}{[1+(2^{-k}r)^2]^{\frac{r+\tau}{2}}}\\
& = C \| \mu \|_r \sum_{k=0}^\infty 2^{-k(r-n)} \int_0^\infty \frac{t^{n-1}dt}{(1+t^2)^{\frac{r+\tau}{2}}}\\
& = C(\| \mu \|_r, r, n, \tau).
\end{align*}
The general case, with $m$ variables, is treated similarly: instead of ordinary sums and integrals one encounters multiple
sums and integrals; we leave details to the reader. $\Box$

\begin{rem}
In the above theorem it is not possible to replace integration over $\mathbb H_{w_j}$ with integration over $\mathbb H$, i.e.
the global variant of this theorem is not true. In fact, a counterexample is obtained by taking $m=1$, $n = 1$, $r =2$,
$\tau = 1$ and $\mu = \sum_{k\geq 1} 2^{2k} \delta_{z_k}$, where $z_k = (0, 2^k)$.
\end{rem}

\section{Multipliers between spaces of harmonic functions on the unit ball}


The goal of this section is to extend  our previous results on multipliers, see \cite{AS3}, to more general harmonic function spaces, involving derivatives. We restrict ourselves to the three theorems below, though other results from our previous work can be generalized similarly. Let us recall some standard notation and facts on spherical harmonics, see \cite{StW} for a
detailed exposition.

Let $Y^{(k)}_j$ be the spherical harmonics of order $k$, $j \leq 1 \leq d_k$, on $\mathbb S$.
The spherical harmonics $Y^{(k)}_j$, ($k \geq 0$, $1 \leq j \leq d_k$), form an orthonormal basis of $L^2(\mathbb S, dx')$. Every $f \in h(\mathbb B)$ has an expansion
$$f(x) = f(rx') = \sum_{k=0}^\infty r^k b_k\cdot Y^k(x'),$$
where $b_k = (b_k^1, \ldots, b_k^{d_k})$, $Y^k = (Y_1^{(k)}, \ldots, Y_{d_k}^{(k)})$ and $b_k\cdot Y^k$ is interpreted in the scalar product sense: $b_k \cdot Y^k = \sum_{j=1}^{d_k} b_k^j Y_j^{(k)}$. We often write, to stress dependence on a function $f \in h(\mathbb B)$, $b_k = b_k(f)$ and $b_k^j = b_k^j(f)$, in fact we have linear functionals $b_k^j$,
$k \geq 0, 1 \leq j \leq d_k$, on the space $h(\mathbb B)$.

We denote the Poisson kernel for the unit ball by $P(x, y')$, it is given by
\begin{align*}
P(x, y') = P_{y'}(x) & = \sum_{k=0}^\infty r^k \sum_{j=1}^{d_k} Y^{(k)}_j(y') Y^{(k)}_j(x') \\
& = \frac{1}{n\omega_n} \frac{1-|x|^2}{|x-y'|^n}, \qquad x = rx' \in \mathbb B, \quad y' \in \mathbb S.
\end{align*}
The Bergman kernel for the harmonic Bergman space $A^p_m$, $m>-1$ is the following function
\begin{equation*}
Q_m(x, y) = 2 \sum_{k=0}^\infty \frac{\Gamma(m + 1 + k + n/2)}{\Gamma(m + 1) \Gamma(k + n/2)}
r^k \rho^k Z_{x'}^{(k)}(y'), \qquad x = rx', \; y = \rho y' \in \mathbb B,
\end{equation*}
see \cite{AS3} and references therein for estimates of this kernel.

Let us recall some definitions from \cite{AS1}.

\begin{df}
For a double indexed sequence of complex numbers
$$c = \{ c_k^j : k \geq 0, 1 \leq j \leq d_k \}$$
and a harmonic function $f(rx') = \sum_{k=0}^\infty r^k \sum_{j=1}^{d_k} b_k^j(f) Y^{(k)}_j(x')$ we define
$$(c \ast f) (rx') = \sum_{k=0}^\infty \sum_{j=1}^{d_k} r^k c_k^j b_k^j(f) Y^{(k)}_j(x'), \qquad rx' \in \mathbb B,$$
if the series converges in $\mathbb B$. Similarly we define convolution of $f, g \in h(\mathbb B)$ by
$$(f \ast g)(rx') = \sum_{k=0}^\infty \sum_{j=1}^{d_k} r^k b_k^j(f)b_k^j(g) Y_j^{(k)}(x'), \qquad rx'\in \mathbb B,$$
it is easily seen that $f \ast g$ is defined and harmonic in $\mathbb B$.
\end{df}

\begin{df}
For $t > 0$ and a harmonic function $f(x) = \sum_{k=0}^\infty r^k b_k(f)\cdot Y^k(x')$ on $\mathbb B$ we define a
fractional derivative of order $t$ of $f$ by the following formula:
$$(\Lambda_t f)(x) = \sum_{k=0}^\infty r^k \frac{\Gamma(k+n/2 + t)}{\Gamma(k+n/2)\Gamma(t)}b_k(f)\cdot Y^k(x'),
\qquad x = rx' \in \mathbb B.$$
\end{df}

Clearly, for $f \in h(\mathbb B)$ and $t>0$ the function $\Lambda_t f$ is also harmonic in $\mathbb B$.

\begin{df}
Let $X$ and $Y$ be subspaces of $h(\mathbb B)$. We say that a double indexed sequence $c$ is a multiplier from $X$ to $Y$ if $c \ast f \in Y$ for every $f \in X$. The vector space of all multipliers from $X$ to $Y$ is denoted by $M_H(X, Y)$.
\end{df}

We associate to such a sequence $c$ a harmonic function
\begin{equation}\label{gc}
g_c(x) = g(x) = \sum_{k\geq 0} r^k \sum_{j=1}^{d_k} c_k^j Y^{(k)}_j(x'), \qquad x = rx' \in \mathbb B,
\end{equation}
and express our conditions in terms of fractional derivatives of $g_c$.



\begin{lem}[\cite{Al1}]\label{subh}
If $f: \Omega \rightarrow \mathbb R$ is harmonic in $\Omega \subset \mathbb R^n$ and if $N \in \mathbb N$, then
$|\nabla^N f|^p$ is subharmonic for $p \geq \frac{n}{n+N}$.
\end{lem}
In particular, $|\nabla^N f|$ is subharmonic and hence $M_1(\nabla^N f, r)$ is increasing for any
$f\in h(\mathbb B)$.

The following three theorems have derivative free counterparts, see \cite{AS3}.

\begin{thm}\label{muldb}
Let $1<s<\infty$, $\alpha, \beta > 0$, $N \in \mathbb N$, $m > \alpha - 1$ and $0<p\leq 1$.
Then $c \in M_H(D_NB^{1,p}_\alpha, H^s_\beta)$ if and only if the function $g = g_c$ satisfies the following condition
\begin{equation}\label{eqmul3}
L_s(g) = \sup_{0\leq\rho < 1} \sup_{y' \in \mathbb S} (1-\rho)^{m +1+N + \beta -\alpha}
\left( \int_{\mathbb S} |\Lambda_{m+1}(g \ast P_{x'})(\rho y')|^s dx' \right)^{1/s} < \infty.
\end{equation}
\end{thm}

{\it Proof.} In proving sufficiency of the condition (\ref{eqmul3}) we follow closely arguments presented in the proof of
Theorem 4 from \cite{AS1}. Namely, let us assume $L_s(g) < \infty$, take $f \in D_NB^{1,p}_\alpha$ and set $h = c \ast f$.
Since $\nabla^N h = c \ast \nabla^N f$, Lemma 6 from \cite{AS1} gives
\begin{equation}\label{tool}
\nabla^N h(r^2 x') = 2 \int_0^1 \int_{\mathbb S} \Lambda_{m+1} (g \ast P_{\xi})(rR x') \nabla^N f(r R \xi)
(1-R^2)^m R^{n-1} d\xi dR
\end{equation}
and this allows us to obtain the following estimate:
\begin{align*}
M_s(\nabla^N h, r^2) & \leq 2 \int_0^1 (1-R^2)^m R^{n-1} \\
&\phantom{\leq 2} \left\| \int_{\mathbb S} \Lambda_{m+1} (g \ast P_\xi)(rRx') \nabla^N f(rR\xi) d \xi \right\|_{L^s(\mathbb S, dx')} dR \notag\\
& \leq 2 \int_0^1 (1-R^2)^m R^{n-1} M_1(\nabla^N f, rR) \sup_{\xi \in \mathbb S} \|
\Lambda_{m+1} (g \ast P_\xi)(rRx') \|_{L^s} dR \notag \\
& \leq C L_s(g) \int_0^1 (1-R)^m M_1(\nabla^N f, rR) (1-rR)^{\alpha - \beta - m - 1 - N} dR \notag \\
& \leq C L_s(g) \int_0^1 M_1(\nabla^N f, rR) (1-rR)^{\alpha - \beta - N - 1} dR \notag \\
& \leq C L_s(g) \int_0^1 M_1(\nabla^N f, rR) \frac{(1-R)^\alpha}{(1-rR)^{\beta + N + 1}} dR.
\end{align*}
Note that $M_1(\nabla^N f, rR)$ is increasing in $0\leq R < 1$, therefore we can combine Lemma 3 from \cite{AS1} and the
above estimate to obtain, for $1/2 \leq r < 1$:
\begin{align*}
M_s^p(\nabla^N h, r^2) & \leq C L_s^p(g) \int_0^1 M_1^p(\nabla^N f, rR)
\frac{(1-R)^{\alpha p + p -1}}{(1-rR)^{p\beta + (N+1)p}} dR  \\
& \leq C L_s^p(g) (1-r)^{-p\beta - Np} \int_0^1 M_1^p(\nabla^N f, R) (1-R)^{\alpha p -1} dR\\
& \leq C L_s^p(g) (1-r)^{-p\beta - Np} \| f \|_{D_NB^{1,p}_\alpha}^p \\
\end{align*}
Therefore $M_s(\nabla^N h, r^2) \leq C L_s(g) (1-r)^{-\beta - N} \| f \|_{D_NB^{1,p}_\alpha}$, which implies $M_s(h, r) \leq
C L_s(g) (1-r)^{-\beta}$. Now we prove necessity of condition (\ref{eqmul3}). Let us set $f_y = Q_m(x,y)$ and $F_y(x) =
\nabla^N f_y(x) = \nabla^N_x Q_m(x,y)$, $x, y \in \mathbb B$. Then using estimate
$$\nabla^N_x Q_m(x,y)| \leq C |\rho x - y'|^{-n-N-m}, \qquad x = rx', y = \rho y', \quad x',y' \in \mathbb S$$
we obtain $M_1(F_y, r) \leq C(1-|y|r)^{-m-N-1}$. Hence $\| F_y \|_{B^{1,p}_\alpha} \leq C (1-|y|)^{\alpha - m - 1 - N}$ which means $\| f_y \|_{D_NB^{1,p}_\alpha} \leq C (1-|y|)^{\alpha - m -1-N}$. Setting $h_y = M_cf_y$ one obtains, as in Lemma 8 from
\cite{AS1}, the estimate
$$\left( \int_{\mathbb S} |\Lambda_{m+1}(g \ast P_{x'})(\rho y')|^s dx' \right)^{1/s} \leq
(1-|y|)^{-\beta} \| h_y \|_{H^s_\beta}.$$
Since, by continuity of $M_c$, $\| h_y \|_{H^s_\beta} \leq C \| f_y \|_{DB^{1,p}_\alpha}$ the proof is completed by combining
the above estimates. $\Box$

Since $D_NA^p_\alpha = D_NB^{p,p}_{\frac{\alpha + 1}{p}}$, taking $p = 1$ we obtain the following corollary.

\begin{cor}
Let $1<s<\infty$, $\alpha, \beta > 0$, $N \in \mathbb N$ amd $m > \alpha - 1$.
Then $c \in M_H(D_NA^1_\alpha, H^s_\beta)$ if and only if the function $g = g_c$ satisfies the following condition
\begin{equation}\label{eqmul4}
K_s(g) = \sup_{0\leq\rho < 1} \sup_{y' \in \mathbb S} (1-\rho)^{m + N + \beta -\alpha}
\left( \int_{\mathbb S} |\Lambda_{m+1}(g \ast P_{x'})(\rho y')|^s dx' \right)^{1/s} < \infty.
\end{equation}
\end{cor}

Analogously to the proof of Theorem \ref{muldb}, one can modify proofs presented in \cite{AS1} and \cite{AS3} to obtain
the following two theorems.

\begin{thm}\label{bpbp}
Let $1 \leq p \leq q \leq \infty$, $1 \leq s \leq \infty$, $N \in \mathbb N$ and $m > \alpha - 1$. Then for a double indexed sequence $c = \{ c_k^j : k \geq 0, 1 \leq j \leq d_k \}$ the following conditions are equivalent:

1. $c \in M_H(D_NB^{1,p}_\alpha, B^{s,q}_\beta)$.

2. The function $g(x) = \sum_{k\geq 0} r^k \sum_{j=1}^{d_k} c_k^j Y^{(k)}_j(x')$ is harmonic in $\mathbb B$ and
satisfies the following condition
\begin{equation}\label{ngs}
N_s(g) = \sup_{0\leq\rho < 1} \sup_{y' \in \mathbb S} (1-\rho)^{\beta - \alpha + m + N + 1}
\left( \int_{\mathbb S} |\Lambda_{m+1}(g \ast P_{x'})(\rho y')|^s dx' \right)^{1/s} < \infty.
\end{equation}
\end{thm}

\begin{thm}
Let $0<p\leq 1 \leq q \leq \infty$, $N \in \mathbb N$ and $m > \alpha -1$. Then for a double indexed sequence
$c = \{ c_k^j : k \geq 0, 1 \leq j \leq d_k \}$ the following conditions are equivalent:

1. $c \in M_H(D_NB^{1,p}_\alpha, F^{q,1}_\beta)$.

2. The function $g(x) = \sum_{k\geq 0} r^k \sum_{j=1}^{d_k} c_k^j Y^{(k)}_j(x')$ is harmonic in $\mathbb B$ and
satisfies the following condition
\begin{equation}\label{mult2}
N_1(g) = \sup_{0\leq\rho < 1} \sup_{y' \in \mathbb S} (1-\rho)^{\beta - \alpha + m + N + 1}
\int_{\mathbb S} |\Lambda_{m+1}(g \ast P_{x'})(\rho y')| dx' < \infty.
\end{equation}
\end{thm}

\end{document}